\documentclass[12pt,twoside]{article}
\usepackage{amssymb}
\usepackage{amsmath}

\date{}
\usepackage{graphicx}

\begin{document}

\centerline{\bf }

\centerline{}

\noindent{\bf\Large Fuzzy semiprime subsets of\\ordered groupoids 
(groupoids)}\bigskip

\bigskip

\noindent{\bf Niovi Kehayopulu$^{\bf 1}$ and Michael Tsingelis$^{\bf 
2}$}\bigskip

\bigskip

\noindent$^1$ {\small\it University of Athens, Department of 
Mathematics, 15784 Panepistimiopolis, Greece}\smallskip

\noindent$^2$ {\small\it Hellenic Open University, School of Science 
and Technology, Studies in Natural Sciences, Greece}\medskip

\bigskip

\begin{center}
{\bf Abstract}
\end{center}{\small
\noindent A fuzzy subset $f$ of an ordered
semigroup (or semigroup) $S$ is called fuzzy semiprime if $f(x)\ge
f(x^2)$ for every $x\in S$ (Definition 1). Following the terminology 
of
semiprime subsets of ordered semigroups (semigroups), the
terminology of ideal elements of $poe$-semigroups (: ordered
semigroups possessing a greatest element), and the terminology of 
ordered semigroups, in general, a fuzzy subset $f$ of
an ordered semigroups (semigroup) should be called fuzzy semiprime
if
for every fuzzy subset $g$ of $S$ such that $g^2:=g\circ g\preceq f$, 
we have
$g\preceq f$ (Definition 2). And this is because if $S$ is a 
semigroup or ordered semigroup, then the set of all fuzzy subsets of 
$S$ is a semigroup (ordered semigroup) as well. What is the relation
between these two definitions? that is between the usual definition 
(Definition 1) we always use and the definition we give in the 
present paper (Definition 2) saying that that definition should 
actually be the correct one? The present paper gives the related
answer.

\bigskip

\smallskip

}
\section{Introduction}
\noindent An ordered groupoid (: $po$-groupoid), denoted by 
$(S,.,\le)$, is an
ordered
set $(S,\le)$ endowed with a multiplication "." which is
compatible with the ordering (that is, $a\le b$ implies $ac\le bc$
and $ca\le cb$ for every $c\in S$). If this multiplication is 
associative, then $S$ is called an ordered semigroup (: 
$po$-semigroup) [1, 2]. A $poe$-semigroup is an ordered
semigroup having a greatest element usually denoted by $"e"$ ($e\ge
a$ for all $a\in S$) (cf. for example [3]). Following L. Zadeh [4], 
the founder of fuzzy sets, if $S$ is an ordered groupoid (or 
groupoid), a fuzzy
subset of $S$ (or a fuzzy set in $S$) is a mapping $f$ of $S$ into
the closed interval $[0,1]$ of real numbers. For a nonempty subset 
$A$ of an ordered groupoid (or groupoid) $S$, the characteristic 
function $f_A$ is the fuzzy subset on $S$ defined 
by\begin{center}$f_A : S\rightarrow \{0,1\} \mid x\rightarrow \left\{ 
\begin{array}{l}
1\,\,\,\,\,if\,\,\,\,\,x \in A\\
0\,\,\,\,\,if\,\,\,\,\,x \notin A
\end{array} \right.$ \end{center}If $(S,.)$ is a groupoid, $f,g$ 
fuzzy subsets of $S$ and $x\in S$, we define$$(f\circ g)(x):=\left\{ 
\begin{array}{l}
\mathop {\sup }\limits_{y,z \in S,\;yz = x} [\min \{ f(y),g(z)\}] 
\,\,\mbox { if there exist }y,z \in S \mbox { such that } x = yz\,\\
0\,\,\,\mbox { if there are no } y,z \in S \mbox { such that } x = 
yz.
\end{array} \right.\,\,\,$$
If $S$ is an ordered groupoid (or groupoid), $x\in S$ and $\lambda\in 
[0,1]$, the mapping\begin{center}$x_{\lambda} : S \rightarrow [0,1] 
\mid y \rightarrow \left\{ \begin{array}{l}
\lambda \,\,\,\,\mbox { if }\,\,\,\,y = x\\
0\,\,\,\,\mbox { if }\,\,\,\,y\ne x
\end{array} \right.$\end{center} is called a {\it fuzzy point} of 
$S$. For short, we write $x_{\lambda}^2$ instead of 
$(x^2)_{\lambda}$.\\
If $S$ is an ordered groupoid, then for an element $a$ of $S$, we 
define
$$A_a=\{(x,y)\in S\times S \mid a\le xy\}.$$For two fuzzy subsets $f, 
g$ of $S$, we define the multiplication $f\circ g$ as:
$$\begin{array}{l}
 (f\circ g)(a):= \left\{ \begin{array}{l}
 \mathop {\bigvee }\limits_{(x,y) \in A_a }  \min \{f(x),g(y)\}
 \,\,\,\,\mbox { if }\,\,A_a  \ne \emptyset  \\
 \,\,\,\,\,\,0\,\,\,\,\,\,\,\,\,\,\,\,\,\,\,\,\,\,\,\,\,\,\,\,\,\,\,\,\,\,\,\,\,\,
 \,\,\,\,\,\,\,\,\,\,\;\;\;\;\;\;\mbox { if }\,\,A_a \, = \emptyset, 
\\
 \end{array} \right. \\
 \, \\
 \end{array}$$and the order relation as follows:$$f\preceq g \mbox { 
if and
 only if } f(x)\le g(x) \mbox { for all } x\in S.$$ If $f, g$ are
 fuzzy subsets of $S$ such that $f\preceq g$ then, for every fuzzy
 subset $h$ of $S$, we have $f\circ h\preceq g\circ h$ and $h\circ
 f\preceq h\circ g$. If $S$ is an ordered semigroup, then the
 multiplication of fuzzy subsets of $S$ is associative, so the set
of all
 fuzzy subsets of $S$ with the multiplication and the order above is
 an ordered semigroup, in particular, a $poe$-semigroup [5]. If $S$ 
is an ordered groupoid, then the set of all fuzzy subsets of $S$ is 
a $poe$-groupoid [5]. Just for an information, if $S$ is an ordered 
groupoid (resp. ordered semigroup), then the $poe$-groupoid (resp. 
$poe$-semigroup) of all fuzzy subsets of $S$ has a zero element and 
$S$ is embedded in the set of all fuzzy subsets of $S$ [5].

According to Clifford and Preston [6; p. 121], a subset $T$ of a 
semigroup $S$ is called semiprime if for
every $a\in S$ such that $a^2\in T$, we have $a\in T$. Semiprime 
ideals play an important role in studying the structure of 
semigroups. As an example, a semigroup $S$ is left regular (resp. 
right regular) if and only if every left (resp. right) ideal of $S$ 
is semiprime. A semigroup $S$ is intra-regular if and only if every 
ideal (that is, two-sided ideal) of $S$ is semiprime. These, in turn, 
are equivalent to saying that a semigroup $S$ is left regular if and 
only it is a union (or disjoint union) of left simple subsemigroups 
of $S$ (the right analogue also holds). Every left and every right 
ideal of $S$ is semiprime if and only $S$ is union of groups (or 
disjoint groups) (which means that every left and every right ideal 
of $S$ is semiprime). A semigroup $S$ is intra-regular (which means 
that the principal ideals of $S$ constitute a semilattice $Y$ under 
intersection) if and only if it is a union is simple semigroups. The
semiprime subsets of ordered semigroups (groupoids) have been defined 
in [7] in the same way. According to [7], a subset $T$ of an ordered 
semigroup $S$ is called semiprime if for every $A\subseteq T$ such 
that $A^2\subseteq T$, we have $A\subseteq T$ (which actually is the 
same with that one given by Clifford and Preston as the two 
definitions are equivalent). In a series of papers the authors of the 
present paper have shown that, exactly as in semigroups, semiprime 
subsets of ordered semigroups play an important role in studying the 
structure of ordered semigroups.

A fuzzy subset $f$ of a semigroup $S$ is called {\it semiprime} if 
$f(x)\ge f(x^2)$ for every $x\in S$. This concept has been first 
introduced by N. Kuroki in [8], as he was the first who observed and 
showed in [8] that a nonempty subset $A$ of a semigroup $S$ is 
semiprime if and only if its characteristic function $f_A$ is fuzzy 
semiprime. Kehayopulu and Tsingelis were the first who studied fuzzy 
ordered groupoids [9]. Following Kuroki, they kept the same 
definition of semiprime subset of an ordered groupoid as a fuzzy 
subset $f$ of $S$ satisfying $f(x)\ge f(x^2)$ for every $x\in S$ [9]. 
Many papers on semigroups and ordered semigroups appeared adapting 
this definition as the definition of semiprime fuzzy subsets both for 
semigroups and ordered semigroups. It might be also noted that a 
fuzzy subset $f$ of a groupoid $S$ is semiprime if and only if for 
every $x\in S$ and every $\lambda\in [0,1]$ such that 
$x_{\lambda}\circ x_{\lambda}\le f$ implies $x_{\lambda}\le f$ [10].

On the other hand, an element $t$ of a $poe$-groupoid (or 
$poe$-semigroup) $S$ is called semiprime if
for every $a\in S$ such that $a^2\le t$, we have $a\le t$ [3]. And 
the same definition of semiprime elements is the usual definition for 
ordered semigroups in general. As this is the case for ordered 
semigroups, in addition, since the fuzzy subsets of ordered 
semigroups form an ordered semigroup, one should expect that in the 
theory of fuzzy ordered semigroups (or semigroups) the fuzzy 
semiprime subset should be defined in a similar way. That
is, if $S$ is an ordered semigroup (or semigroup), then a fuzzy
subset $f$ of $S$ should be called fuzzy semiprime if for any fuzzy
subset $g$ of $S$ such that $g\circ g\preceq f$, we have $g\preceq
f$. However in the existing bibliography, for an ordered semigroup
$S$, a fuzzy subset $f$ of $S$ is called fuzzy semiprime if $f(x)\ge
f(x^2)$ for every $x\in S$ (cf. for example [8--12]) and this is the 
usual definition the authors always use. It is natural to ask what is 
the
relation
between these two definitions. The present paper gives the related 
answer. Here we prove that if
a fuzzy subset $f$ of an ordered groupoid (semigroup) $S$ is 
semiprime (in the
usual sense), then for any fuzzy subset $g$ of $S$ such that $g\circ
g\preceq f$, we have $g\preceq f$, and that the converse statement
does not hold in general. However, for the ordered semigroups
satisfying the condition $$x\le yz \Longrightarrow \min \{f(y^2),
f(z^2)\}\le f(x),$$ the two definitions are equivalent.
\section{Main results}
\noindent{\bf Proposition 1.} {\it Let $(S,.,\le)$ be an ordered
groupoid and f a fuzzy subset of S. Then$$f(x)\le (f\circ f)(x^2)
\mbox { for every } x\in S.$$}{\bf Proof.} Let $x\in S$. Since
$(x,x)\in A_{x^2}$, we have $A_{x^2}\not=\emptyset$ and
\begin{eqnarray*}(f\circ f)(x^2)&=&\mathop  \bigvee \limits_{(u,v)
\in A_{x^2} } \min
\{f(u),g(v)\}\\&\geqslant&\min\{f(x),f(x)\}\\&=&f(x),
\end{eqnarray*}so $f(x)\le (f\circ f)(x^2)$. $\hfill\Box$\medskip

\noindent{\bf Proposition 2.} {\it Let $(S,.,\le)$ be an ordered
groupoid and f, g fuzzy subsets of S such that $g\circ g\preceq f$.
Then$$g(x)\le f(x^2) \mbox { for every } x\in S.$$}{\bf Proof.} Let
$x\in S$. Since $g$ is a fuzzy subset of $S$, by Proposition 1, we
have $g(x)\le (g\circ g)(x^2)$. Since $g\circ g\preceq f$, we have
$(g\circ g)(x^2)\le f(x^2)$. Thus we have $g(x)\le f(x^2)$.
$\hfill\Box$\medskip

\noindent{\bf Definition 3.} [9] If $S$ is an ordered groupoid, a 
fuzzy
subset $f$ of $S$ is called {\it fuzzy semiprime} if $f(x)\ge
f(x^2)$
for every $x\in S$.\medskip

\noindent{\bf Theorem 4.} {\it Let S be an ordered groupoid and f a
fuzzy subset of S. We consider the following statements:

$(1)$ f is fuzzy semiprime.

$(2)$ If g is a fuzzy subset of S such that $g\circ g\preceq f$,
then
$g\preceq f$.\\Then $(1)\Rightarrow (2)$. The implication
$(2)\Rightarrow (1)$ does not hold in general. }\medskip

\noindent{\bf Proof.} $(1)\Longrightarrow (2)$. Let $g$ be a fuzzy
subset of $S$ such that $g\circ g\preceq f$ and $x\in S$. By
Proposition 2, we have $g(x)\le f(x^2)$. Since $f$ is fuzzy
semiprime, we have $f(x)\ge f(x^2)$, Then we have $g(x)\le f(x)$ and
(2) holds.\medskip

\noindent Condition (2) does not always imply (1). In fact:\medskip

The set $S=\{n\in N \mid n\ge 2\}=\{2,3,4,.....\}$ of natural
numbers
with the usual multiplication and the usual order is an ordered 
groupoid (in particular, it is an ordered semigroup). Let
$f$ be the fuzzy subset of $S$ defined by:
$$\,f\,\,:\,\,\,(S,., \le )\,\, \to \,\,[0,1] \mid x \to \left\{
\begin{array}{l}
0\,\,\,\,\,{\rm{if}}\,\,\,\,\,x = 2\\
1\,\,\,\,\,{\rm{if}}\,\,\,\,\,x > 2.
\end{array} \right.$${\it Condition $(2)$ is satisfied.} Indeed: Let
$g$ be a fuzzy subset of $S$ such that $g\circ g\preceq f$ and let
$x\in S$. \medskip

(I) Let $x=2$.

We consider the set $A_2=\{(m,n)\in S\times S \mid 2\le mn\}$. Since
$(2,2)\in A_2$, we have $A_2\not=\emptyset$ and
\begin{eqnarray*}(g\circ g)(2)&=&\mathop  \bigvee \limits_{(m,n) \in
A_{2} } \min \{g(m),g(n)\}\\&\geqslant&\min\{g(2),g(2)\}\\&=&g(2),
\end{eqnarray*}Since $g\circ g\preceq f$, we have $(g\circ g)(2)\le
f(2)$, so we have $g(2)\le f(2)$.\medskip

(II) Let $x>2$. Then $f(x)=1$. On the other hand, since $g$ is a
fuzzy subset of $S$, we have $g(x)\le 1$. Thus we have $g(x)\le
f(x)$.\medskip

By (I) and (II), we have $g(x)\le f(x)$. Since this holds for any
$x\in S$, we have $g\preceq f$.\medskip

\noindent {\it Condition $(1)$ does not hold.} In fact, we have
$f(2)=0$ and $f(2^2)=f(4)=1$, so $f(2)\not\ge f(2^2)$.
$\hfill\Box$\medskip

It is natural to ask under what conditions the implication
$(2)\Rightarrow (1)$ is satisfied. The next theorem gives a related
answer. \medskip

\noindent{\bf Theorem 5.} {\it Let S be an ordered groupoid and f a
fuzzy subset of S such that}

(a) {\it $a\le xy$ $\Longrightarrow$ $\min\{f(x^2),f(y^2)\}\le f(a)$
$(a,x,y\in S)$ and}

(b) {\it if g is a fuzzy subset of S such that $g\circ g\preceq f$,
then $g\preceq f$.\\Then f is fuzzy semiprime}.\medskip

\noindent{\bf Proof.} Let $x\in S$. We consider the fuzzy subset $g$
of $S$ defined by:$$g : (S,.,\le) \rightarrow [0,1] \mid x
\rightarrow f(x^2).$$We have $g\circ g\preceq f$. In fact: Let $a\in
S$. If $A_a=\emptyset$, then $(g\circ g)(a)=0\le f(a)$. If
$A_a\not=\emptyset$, then $$(g\circ g)(a)=\mathop  \bigvee
\limits_{(x,y) \in A_{2} } \min \{g(x),g(y)\}.$$We have
$$\min\{g(x),g(y)\}\le f(a) \mbox { for every } (x,y)\in A_a.$$
Indeed, if $(x,y)\in A_a$, then $a\le xy$ and, by (a),
$\min\{f(x^2),f(y^2)\}\le f(a)$, that is $\min\{g(x),g(y)\}\le
f(a)$.
Therefore we have $(g\circ g)(a)\le f(a)$. This is for every $a\in
S$, thus we obtain $g\circ g\preceq f$. By condition (a), we get
$g\preceq f$, then $f(x)\ge g(x)=f(x^2)$. This holds for every $x\in
S$, so $S$ is fuzzy semiprime$.\hfill\Box$\medskip

\noindent{\bf Remark 6.} Fuzzy semiprime subsets of ordered
semigroups do not satisfy the condition (a) of Theorem 5 in general.
In fact: Let $S=[0,1]$ be the ordered semigroup of real numbers with
the usual multiplication and the usual order of reals and $f$ the 
fuzzy subset on
$S$ defined by$$f : (S,.,\le) \rightarrow [0,1] \mid x \rightarrow
x$$ (that is, the identity mapping on $S$). Then {\it $f$ is fuzzy
semiprime}. Indeed: If $x\in S$, then $x<1$. Since $x\ge 0$, we have
$x^2\le x$, so $f(x)\ge f(x^2)$.

\noindent{\it $S$ does not satisfy the condition} (a). In fact:
$$\frac{1}{{10}}\le \frac{1}{{6}}=\frac{1}{{2}}.\frac{1}{{3}}$$but
\begin{eqnarray*}\min\{f((\frac{1}{{2}})^2),f((\frac{1}{{3}})^2)\}
&=&\min\{f(\frac{1}{{4}}),f(\frac{1}{{9}})\}=\min\{\frac{1}{{4}},
\frac{1}{{9}}\}\\&=&\frac{1}{{9}}\not\le\frac{1}{{10}}=f(\frac{1}{{10}}).
\end{eqnarray*}\section{Conclusion}
As as conclusion, let us give the two definitions below: The first
one is the definition in the existing bibliography we always use. The 
second is
similar with the definition of semiprime subsets (or ideal elements)
of ordered groupoid. The definition which should actually
be.\smallskip

In the following, $S$ is an ordered groupoid (or groupoid) and
$g^2:=g\circ g$.\medskip

\noindent{\bf Definition 1.} A fuzzy subset $f$ of $S$ is called
{\it fuzzy semiprime} if $$f(x)\ge f(x^2)$$for every $x\in 
S$.\medskip

\noindent{\bf Definition 2.} A fuzzy subset $f$ of $S$ is called
{\it fuzzy semiprime} if\smallskip

For any fuzzy subset $g$ of $S$ such that $g^2\preceq f$, we have
$g\preceq f$. \medskip

\noindent Then Definition 1 implies Definition 2, but Definition 2 
does not
imply Definition 1 in general. In particular, if the fuzzy subset
$f$ of $S$
has the property$$a\le xy \Longrightarrow \min \{f(x^2), f(y^2)\}\le
f(a),$$then the two definitions are equivalent.\bigskip

\noindent This paper has been submitted in International Journal of 
Mathematics and Mathematical Sciences on January 2, 2014

\end{document}